\documentclass[11pt,a4paper]{article}

\usepackage{epsfig,amsfonts,amsgen,
amsmath,amstext,amsbsy,amsopn,amsthm,amssymb,mathrsfs, epstopdf}
\usepackage{kantlipsum,multirow}
\usepackage{tikz,chemfig}
\allowdisplaybreaks

\usetikzlibrary{backgrounds}
\usetikzlibrary{arrows}
\usetikzlibrary{shapes,shapes.geometric,shapes.misc}
\pgfdeclarelayer{edgelayer}
\pgfdeclarelayer{nodelayer}
\pgfsetlayers{background,edgelayer,nodelayer,main}
\tikzstyle{none}=[inner sep=0mm]
\tikzstyle{every loop}=[]

\tikzstyle{dotted}=[dash pattern=on \pgflinewidth off 2pt]
\tikzstyle{dashed}=[dash pattern=on 3pt off 3pt]

\newcommand \tikzp[2]
{
\begin{center}
\begin{tikzpicture}[scale=#1]
#2
\end{tikzpicture}
\end{center}
}

\tikzstyle{new style 0}=[fill=black, draw=black, shape=circle]
\tikzstyle{red style 1}=[fill=red, draw=black, shape=circle]
\tikzstyle{blue style 2}=[fill=blue, draw=black, shape=circle]
\tikzstyle{white style 4}=[fill=white, draw=black, shape=circle]
\tikzstyle{bklack style 5}=[fill=black, draw=black, shape=rectangle]
\tikzstyle{red style 3}=[fill=red, draw=black, shape=rectangle]
\tikzstyle{yellow style 7}=[fill=yellow, draw=black, shape=rectangle]
\tikzstyle{new style 8}=[fill={rgb,255: red,0; green,132; blue,0}, draw={rgb,255: red,0; green,131; blue,0}, shape=circle]

\tikzstyle{new edge style 0}=[-]
\tikzstyle{new edge style 1}=[-, draw=red]
\tikzstyle{new edge style 2}=[-, draw=blue]
\tikzstyle{new edge style 3}=[-, draw={rgb,255: red,0; green,156; blue,0}]

\tikzstyle{cblue}=[circle, draw, thin,fill=blue!20, scale=0.5]

\newtheorem{theorem}{Theorem} 

\newtheorem{lemma}{Lemma}
\newtheorem{false statement}{False statement}

\theoremstyle{definition}

\newtheorem{claim}{Claim}

\newtheorem{corollary}[claim]{Corollary}

\def \proof {\noindent {\it Proof}. }
\def \proofend {\hfill $\Box$}

\def \S {{\mathbb S}}
\def \D {{\mathbb D}}

\def \oar {\overrightarrow}

\newcommand \red[1]{\textcolor{red}{#1}}

\def \T {\mathcal{T}}

\def \NST {\mathcal{NST}}

\newcommand \equ[2]
{
\begin{equation}\label{#1}
#2
\end{equation}
}

\newcounter{countclaim}
\def\inclaim{\addtocounter{countclaim}{1}
{\noindent {\bf Claim \thecountclaim}: }}

\newcommand {\rebibitem}[1]
{\bibitem{#1} 
}
\begin{document}

\title
{Express the number of spanning trees in term of degrees\thanks{The work was supported by the
National Natural Science Foundation of China
(No. 11701401) and the Scientific Research
Fund of Hunan Provincial Education Department of China (No. 18A432).}
}

\date{}

\def \bg {\hspace{0.3 cm}}

\author{
Fengming Dong\thanks{Corresponding author.
Email: fengming.dong@nie.edu.sg 
and donggraph@163.com.} \\
\small National Institute of Education, Nanyang Technological University,
Singapore
\\
Jun Ge\thanks{Email: mathsgejun@163.com.}\\
\small School of Mathematical Sciences,
Sichuan Normal University, Chengdu, P. R. China\\
Zhangdong Ouyang\thanks{oymath@163.com.}\\
\small Department of Mathematics,
Hunan First Normal University, Changsha, 
P.R. China
}

\maketitle

\begin{abstract}
It is well-known that the number of spanning trees,
denoted by $\tau(G)$,
in a connected multi-graph $G$ can be calculated
by the Matrix-Tree Theorem and 
Tutte's deletion-contraction formula.
In this short note, we find an alternate method
to compute $\tau(G)$ by degrees of vertices.
\end{abstract}

\medskip



\section{Introduction}\label{intro}

In this article, we consider
loopless and undirected multi-graphs.
For a graph $G$, let
$V(G), E(G)$ and $\T(G)$ be the set of vertices,
the set of edges and the set of spanning trees in $G$
respectively, and let $\tau(G)=|\T(G)|$.
For any $u\in V(G)$,
let $E_G(u)$ (or simply $E(u)$) denote the set of edges in $G$
that are incident with $u$,
and let $d_G(u)$ (or simply $d(u)$) be the degree of $u$ in $G$,
i.e., $d_G(u)=|E_G(u)|$.
For any $S\subseteq V(G)$, if $S\ne \emptyset$,
let $G[S]$ be the subgraph of $G$ induced by $S$,
and if $S\ne V$, let $G-S=G[V\setminus S]$.


The study of spanning trees plays an important role in
graph theory.
The number of spanning trees $\tau(G)$
is a key parameter in Tutte polynomials,
and it has a close relation with some other parameters.
Given a multi-graph $G$,
$\tau(G)=0$ if and only if $G$ is disconnected.
When $G$ is connected,
$\tau(G)$ can be computed by some different methods,
such as 
Kirchhoff's Matrix-Tree Theorem~\cite{Kirch1847,Kocay},
Tutte's deletion-contraction formula~\cite{Tutte}, etc.
In some special cases, $\tau(G)$ can be computed
directly by explicit formulas.
The most famous one is Cayley's formula,
i.e., $\tau(K_n)=n^{n-2}$ for complete graphs \cite{Cayley}.
This formula has been extended to
$\tau(K_{n_1,n_2,\cdots,n_k})=n^{k-2}\prod_{i=1}^k (n-n_i)^{n_i-1}$
for any complete $k$-particle
graph $K_{n_1,n_2,\cdots,n_k}$,
where $n=n_1+n_2+\cdots+n_k$
\cite{Austin}.
It is also known that
$\tau(Q_n)=2^{{2^{n}-n-1}}
\prod _{{k=2}}^{n}k^{{{n \choose k}}}
$ 
for the $n$-dimensional hypercube graph $Q_n$ \cite{Harary}.
For the line graph $G=L(H)$ of an arbitrary connected graph $H$,
a relation between $\tau(G)$ and spanning trees in
$H$ was also established~\cite{DongYan}.
More works on $\tau(G)$ can be found in
\cite{Ge,Gong,Lovasz,Moon,Yan}.

In the following is an upper bounds for $\tau(G)$ due to
Thomassen \cite{Thomassen2010}.

\begin{theorem}[\cite{Thomassen2010}]\label{Thomassen}
Let $G=(V,E)$ be a multi-graph and $u$ be any vertex in $G$.
Then
$$\tau (G)\leq \prod_{v\in V-\{u\}}d(v).$$
\end{theorem}

For any multi-graph $G$ and any vertex $u$ in $G$,
let $\NST_u(G)$ be the set of non-spanning subtrees $T$ of
$G$ such that 
$u\in V(T)$ and $G-V(T)$ has no isolated vertices.
In this article, we find the following formula
expressing $\tau(G)$ in terms of degrees.
It shows how far is Thomassen's
upper bound from $\tau(G)$ exactly.

\begin{theorem}\label{main}
For a multi-graph $G=(V,E)$ and a vertex $u$ in $G$,
\equ{main-eq1}
{
\tau(G)=\prod_{v\in V-\{u\}}d(v)-\sum_{T\in \NST_u(G)}\prod_{v\in V-V(T)}d_{G-V(T)}(v).
}
\end{theorem}

Theorem \ref{main} can be proved by 
some different approaches.
In this note, we shall prove Theorem~\ref{identity}
in Section 3
from which Theorem~\ref{main} follows directly.
In Section~\ref{sec2},
we introduce a polynomial $F(G,\omega)$ of a graph $G$
by assigning a variable $y_i$ to each edge $e_i$ in $G$.
This polynomial will be applied in Section~\ref{sec3}
for proving Theorem~\ref{identity} by
a method inspired by Wang algebra \cite{Duffin1959, Wang1934}\footnote{Wang algebra
assumes that $x+x=0,x\cdot x=0$ and $xy=yx$
for any variables $x$ and $y$.}.
In Section~\ref{sec4},
we apply Theorem~\ref{main} to compute $\tau(G)$
for some graphs.


\section{A polynomial $F(G,\omega)$}\label{sec2}

For any positive integer $n$, let $[n]=\{1,2,\cdots,n\}$.
Let  $G=(V,E)$ be a loopless and connected multi-graph
with $V=\{v_i: i\in [n]\}$ and $E=\{e_j: j\in [m]\}$.
Assume that $\omega$ is a weight function on $E(G)$
defined by $\omega(e_j)=y_j$ for each $j\in [m]$,
where $y_1,y_2,\cdots,y_m$ are considered as indeterminates.
Define a polynomial $F(G,\omega)$ as follows:
\equ{eq1}
{
F(G,\omega)
=\prod_{i\in [n]} \sum_{e_j\in E(v_i)}\omega(e_j)
=\prod_{i\in [n]} \sum_{e_j\in E(v_i)}y_j,
\qquad
\mbox{when } V\ne \emptyset;
}
and $F(G, \omega)=1$ when $V=\emptyset$.
Clearly, $F(G,\omega)=0$ whenever $d(v_i)=0$ for some $v_i\in V$.
If $y_i=1$ for all $i\in [m]$, then
$
F(G,\omega)=\prod\limits_{1\le i\le n} d_G(v_i).
$

\def \F {{\mathscr F}}

The expansion of $F(G,\omega)$ can be applied to study
some structures of $G$,
such as 
the minimum edge coverings, maximum matchings,
perfect matchings, and spanning trees,
and hence the edge covering number $\rho(G)$,
the matching number $\nu(G)$ and
the number of spanning trees $\tau(G)$.
Let $\F(G,\omega)$ denote the set of terms in the expansion of $F(G,\omega)$.
Note that each term in $\F(G,\omega)$
is in the form
$y_{i_1}^2y_{i_2}^2\cdots y_{i_r}^2y_{j_1}\cdots y_{j_k}$,
where $k+2r=n$ and
$i_1,i_2,\cdots,i_r,j_1,\cdots, j_k$ are pairwise distinct.
Each term
$y_{i_1}^2y_{i_2}^2\cdots y_{i_r}^2y_{j_1}\cdots y_{j_k}$
in $\F(G,\omega)$
corresponds to an edge cover $\{e_{i_1},e_{i_2},
\cdots,e_{i_r}\}\cup \{e_{j_1},e_{j_2},\cdots, e_{j_k}\}$
of $G$, where
$\{e_{i_1},e_{i_2}, \cdots,e_{i_r}\}$ is a matching of $G$.
In particular, if
$y_{i_1}^2y_{i_2}^2\cdots y_{i_r}^2$ is a term
in $\F(G,\omega)$, 
then $n=2r$ and
it corresponds to a perfect matching
$\{e_{i_1},e_{i_2}, \cdots,e_{i_r}\}$ of $G$.
Thus, $\rho(G)$ is the minimum value of $k+r$
among all terms
$y_{i_1}^2y_{i_2}^2\cdots y_{i_r}^2y_{j_1}\cdots y_{j_k}$
in $\F(G,\omega)$,
and $\nu(G)$ is the maximum value of $r$
among all terms $y_{i_1}^2y_{i_2}^2\cdots y_{i_r}^2y_{j_1}\cdots y_{j_k}$
in $\F(G,\omega)$.

\begin{figure}[htp]
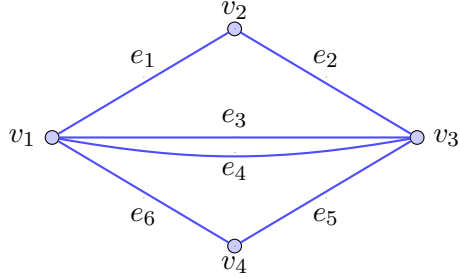

\tikzp{0.8}
{
\foreach \place/\x in {{(-3,0)/1}, {(0,1.8)/2},{(3,0)/3}, {(0,-1.8)/4}}
\node[cblue] (v\x) at \place {};

\filldraw[black] (-3.1,0) circle (0pt) node[anchor=east] {$v_1$};
\filldraw[black] (v2) circle (0pt) node[anchor=south] {$v_2$};
\filldraw[black] (3.1,0) circle (0pt) node[anchor=west] {$v_3$};
\filldraw[black] (v4) circle (0pt) node[anchor=north] {$v_4$};

\filldraw[black] (-1.5,1) circle (0pt) node[anchor=south] {$e_1$};
\filldraw[black] (1.5,1) circle (0pt) node[anchor=south] {$e_2$};
\filldraw[black] (-1.5,-1) circle (0pt) node[anchor=north] {$e_6$};
\filldraw[black] (1.5,-1) circle (0pt) node[anchor=north] {$e_5$};
\filldraw[black] (0,0.01) circle (0pt) node[anchor=south] {$e_3$};
\filldraw[black] (0,-0.25) circle (0pt) node[anchor=north] {$e_4$};

\draw[blue!70,  thick] (v1) -- (v2) -- (v3)--(v4)--(v1)--(v3);
\draw [blue!70,  thick, bend right=10, looseness=1] (v1) to (v3);
}
\caption{A multi-graph}
\label{f1}
\end{figure}

For example, if $G$ is the multi-graph in Figure~\ref{f1}, then
\equ{eq2}
{F(G,\omega)=(y_1+y_3+y_4+y_6)(y_1+y_2)(y_2+y_3+y_4+y_5)(y_5+y_6),
}
and the expansion of $F(G,\omega)$ contains terms $y_1^2y_5^2$
and $y_2^2y_6^2$, which correspond to the two perfect matchings
in $G$: $M_1=\{e_1,e_5\}$ and $M_2=\{e_2,e_6\}$.

In the next section, we shall apply $F(G,\omega)$ to
study $\tau(G)$.

\section{An identity associated with spanning trees }\label{sec3}

In this section, we assume that 
$G=(V,E)$ is a loopless connected multi-graph,
where $V=\{v_i: i\in [n]\}$, $n\ge 2$, 
and $E=\{e_j: j\in [m]\}$.
Let $\omega$ be a weight function on $E$.

We first establish two lemmas which will be applied 
to prove the main result in this section.

Let $\oar{G}$ denote the digraph obtained from $G$ 
by replacing each edge $e_i$ in $G$ by 
two arcs which are incident the same pair of ends of $e_i$
and have opposite directions. 
Assume that  the weight function $\omega$ is extended to 
the arc set $A(\oar{G})$ such that 
$\omega(a)=\omega(e_i)$ for each 
$a\in A(\oar{G})$ if $a$ is obtained from $e_i$
by assigning a direction.

For a digraph $D$ and a vertex $v$ in $D$,
let $id_D(v)$ denote the in-degree of $v$ in $D$.
If $id_D(v)=0$, then $v$ is called a {\it source} of $D$.

Let $\D^*$ denote the family of spanning subdigraphs 
$D$ of $\oar{G}$ 
with $id_D(v_n)=0$ 
and $id_D(v_i)=1$ for each $i\in [n-1]$.

For any subdigraph $D$ of $\oar{G}$, 
let $\omega(D)=\prod\limits_{a\in A(D)} \omega(a)$
if $A(D)\ne \emptyset$ and $\omega(D)=1$ otherwise.

\begin{lemma}\label{le2-1}
Let $G=(V,E)$ be a loopless connected multi-graph,
where $V=\{v_i: i\in [n]\}$, $n\ge 2$ 
and $E=\{e_j: j\in [m]\}$,
and let $\omega$ be a weight function on $E$.
The following holds: 
\equ{le2-1-e2}
{
\prod_{i=1}^{n-1}\sum_{e_j\in E(v_i)} \omega(e_j)
=\sum_{D\in \D^*}\omega(D).
}
\end{lemma}

\proof 
Let $\Pi$ be the set of mappings 
$\pi: [n-1]\rightarrow [m]$ such that 
$e_{\pi(i)}\in E(v_i)$ for each $i\in [n-1]$.
Observe that 
\equ{le2-1-e3}
{
\prod_{i=1}^{n-1} \sum_{e_j\in E(v_i)}\omega(e_j)
=\sum_{\pi\in \Pi} 
\prod_{1\le i\le n-1} \omega(e_{\pi(i)}).
}
For any $\pi\in \Pi$, 
$(e_{\pi(1)},e_{\pi(2)},\cdots,e_{\pi(n-1)})$ is a 
list of $n-1$ edges in $G$, 
where each edge $e_{\pi(i)}$ is incident to $v_i$.
Let $f(\pi)$ denote the spanning subdigraph $D$ of $\oar{G}$ 
that can be obtained 
by converting each edge 
$e_{\pi(i)}$ into an arc with $v_i$ as its head.
Observe that $f(\pi)$ is a digraph in $\D^*$
and,  if $D=f(\pi)$, then  
\equ{le2-1-e4}
{
\prod_{1\le i\le n-1} \omega(e_{\pi(i)})
=\prod_{a\in A(D)}\omega(a)
=\omega(D).
}

It is obvious that $f$  is a bijection 
from $\Pi$ to $\D^*$.
Thus, (\ref{le2-1-e2}) follows 
from (\ref{le2-1-e3}) and (\ref{le2-1-e4})
and the lemma holds.
\proofend

For any $U\subseteq V$ with $U\ne \emptyset$, 
let $\D[U]$ denote the family of subdigraphs 
$D$ of $\oar{G}$ with vertex set $U$ 
and $id_D(v_i)=1$ for each $v_i\in U$.
Note that $\D[V]$ is different from $\D^*$,
although both are spanning subdigraphs of $\oar{G}$.
The following lemma can be proved similarly.

\begin{lemma}\label{le2-2}
Let $G=(V,E)$ be a loopless connected multi-graph,
where $V=\{v_i: i\in [n]\}$, $n\ge 2$ 
and $E=\{e_j: j\in [m]\}$,
and let $\omega$ be a weight function on $E$.
For any $U\subseteq V(G)$ with $U\ne \emptyset$,
\equ{le2-2-e1}
{
F(G[U], \omega)
=\sum_{D\in \D[U]}\omega(D).
}
\end{lemma}

\vspace{0.2 cm}

Recall that $\T(G)$ is the set of spanning trees in $G$.
For any $T\in \T(G)$, let $\tau(T,\omega)=1$
when $|V(G)|=1$, and let
\equ{tree}
{
\tau(T,\omega)
=\prod_{e_i\in E(T)}\omega(e_i),
\qquad \mbox{when }|V(G)|\ge 2. 
}
Now we define another function $\tau(G,\omega)$:
\equ{FT}
{
\tau(G,\omega)=\sum_{T\in \T(G)}\tau(T,\omega).
}
Thus $\tau(G,\omega)=0$ whenever $\T(G)=\emptyset$
(i.e., $G$ is disconnected).
Clearly, when $G$ is connected,
every term in the expansion of $\tau(G,\omega)$
corresponds to a spanning tree in $G$,
and $\tau(G,\omega)=\tau(G)$ whenever 
$\omega(e_j)=1$
for all $j\in [m]$.

Recall that for any $u\in V(G)$, $\NST_u(G)$ denotes
the set of non-spanning subtrees $T$ of $G$ 
such that $u\in V(T)$ and $G-V(T)$ has no isolated vertices.
We are now going to prove the following identity
on $\tau(G,\omega)$
from which Theorem~\ref{main} follows directly.

\begin{theorem}\label{identity}
Let  $G=(V,E)$ be a loopless connected multi-graph,
where $V=\{v_i: i\in [n]\}$, $n\ge 2$ 
and $E=\{e_j: j\in [m]\}$.
Assume that $\omega$ is a weight function on $E$. 
Then,
\equ{id1}
{\prod_{i=1}^{n-1} \sum_{e_j\in E_{G}(v_i)}\omega(e_j)
=\tau(G,\omega)
+\sum_{T_0\in \NST_{v_n}(G)}\tau(T_0,\omega) F(G-V(T_0),\omega).
}
\end{theorem}

\proof
A digraph is called a {\it directed tree} if 
its underlying graph is a tree. 
A directed tree with a unique source 
is called a {\it rooted directed  tree}
and the unique source is its root.
We are now going to establish the following claims.

\inclaim 
For any weakly connected diraph $Q$ with vertices $u_0,u_1,\cdots,u_{k}$,
if $id_Q(u_0)=0$ and $id_Q(u_i)=1$ for all $i\in [k]$,
then $Q$ is a directed rooted tree with root $u_0$.

$Q$ is a directed tree 
as its underlying graph is connected 
and has exactly $k$ edges and $k+1$ vertices. 
Then the claim holds as $u_0$ is the only source in $Q$.

Recall that $\D^*$ is the family of spanning subdigraphs
$D$ of $\oar{G}$ such that $id_D(v_n)=0$
and $id_D(v_i)=1$ for each $i\in [n-1]$.
For any $D\in \D^*$, 
let $D_{v_n}$ denote the component (i.e., a weakly connected 
component)
of $D$ that contains vertex $v_n$.

\inclaim 
For any $D\in \D^*$, $D_{v_n}$  
is a rooted directed tree with root $v_n$.

If $V(D_{v_n})=\{v_n\}$, the claim is trivial.
Now, without loss of generality, assume that
$V(D_{v_n})
=\{v_i:i\in [k]\}\cup \{v_n\}$, where $1\le k\le n-1$.
As $D_{v_n}$ is weakly connected and $|V(D_{v_n})|=k+1$,
we have $|A(D_{v_n})|\ge k$.

It is known that $D$ has exactly $n-1$ arcs and 
$id_D(v_i)=1$ for all $i\in [n-1]$.
Assume that $a_i$ is the arc in $D$ with head $v_i$
for each $i\in [n-1]$.
Thus, $A(D)=\{a_i:i\in [n-1]\}$.
As $V(D_{v_n})=\{v_i:i\in [k]\}\cup \{v_n\}$,
we have $A(D_{v_n})\subseteq \{a_i:i\in [k]\}$.
Since $|A(D_{v_n})|\ge k$,
$A(D_{v_n})=\{a_i:i\in [k]\}$ holds.

Thus, 
$D_{v_n}$ is weakly connected with a source $v_n$ 
and $a_i$ is the only arc in $D_{v_n}$ with head $v_i$
for all $i\in [k]$. 
Claim 2 then follows from Claim 1.

\inclaim For each subtree $T$ of $G$ with $v_n\in V(T)$,
there is exactly one rooted directed tree,
denoted by $\oar{T}$,
with the following properties:
\begin{enumerate}
\item [(i)] $T$ is the underlying graph of $\oar{T}$;
and 

\item [(ii)] $id_{\oar{T}}(v_n)=0$ 
and $id_{\oar{T}}(v_i)=1$ 
for each $v_i\in V(T)\setminus \{v_n\}$.
\end{enumerate}

Claim 3 is obvious, as such a directed tree $\oar{T}$
can only be obtained by assigning directions to
edges in $T$ so that each $v_n-v_i$ path in $T$ 
becomes a directed $v_n-v_i$ path 
(i.e., a path from $v_n$ to $v_i$)
in $\oar{T}$.
Observe that $\omega(T)=\omega(\oar{T})$ for 
each subtree $T$ of $G$.

Recall that $\NST_{v_n}(G)$ is the set of non-spanning 
subtrees $T$ of $G$ such that $v_n\in V(T)$ and 
$G-V(T)$ has no isolated vertices.
Let $\NST_{v_n}(\oar{G})=\{\oar {T}: T\in \NST_{v_n}(G)\}$.

By Claim 2, for each $D\in \D^*$,
if $D$ is not weakly connected, 
then, the unlderlying graph $T$ of 
$D_{v_n}$ is a non-spanning tree.
Furthermore, by the definition of $\D^*$, 
each vertex $v_i$, where $i\in [n-1]$, 
is the head of some arc in $\D^*$
and thus is not isolated in $G-V(T)$,
implying that 
$D_{v_n}=\oar{T}\in \NST_{v_n}(\oar{G})$.

For any $\oar{T}\in \NST_{v_n}(\oar{G})$,
let $\D^*(\oar{T})$ denote the set of 
$D\in \D^*$ such that 
$D_{v_n}$ is the directed tree $\oar{T}$.
Thus, by the definition of $\D[U]$ for $U\subseteq V$, 
for any $T\in \NST_{v_n}(G)$,
\equ
{eq2-7}
{
\D^*(\oar{T}) = \{\oar{T}\cup Q: Q\in \D[V(G)\setminus V(T)]\},
}
where $\oar{T}\cup Q$ denotes the spanning digraph of $\oar{G}$ with arc set $A(\oar{T})\cup A(Q)$.

Let $\D^*_0$ denote the family of $D\in \D^*$ 
such that $D$ is weakly connected.
By Claim 2, $D$ is a rooted directed tree 
for each $D\in \D^*_0$. 
Actually, $\D^*_0=\{\oar{T}: T\in \T(G)\}$.
As $D_{v_n}$ belongs 
to $\NST_{v_n}(\oar{G})$ for each 
$D\in \D^*\setminus \D^*_0$,
by (\ref{eq2-7}), 
\equ
{eq2-5}
{
\D^*\setminus \D^*_0
=\bigcup_{T\in \NST_{v_n}(G)}\D^*(\oar{T})
=\bigcup_{T\in \NST_{v_n}(G)}
\{\oar{T}\cup Q: Q\in \D[V(G)\setminus V(T)]\}.
}

By Lemmas~\ref{le2-1} and~\ref{le2-2} and  (\ref{eq2-5}),
\begin{eqnarray}\label{id1-p}
\prod_{i=1}^{n-1} \sum_{e_j\in E_G(v_i)}\omega(e_j)
&=& \sum_{D\in \D^*_0}\omega(D)
+\sum_{D\in \D^*\setminus \D^*_0}\omega(D)
\nonumber \\
&=&\sum_{T\in \T(G)}\omega(\oar{T})
+\sum_{T\in \NST_{v_n}(G)}
\sum_{Q\in \D[V(G)\setminus V(T)]}
\omega(\oar{T})\omega(Q)
\nonumber \\
&=&\sum_{T\in \T(G)}\omega(T)
+\sum_{T\in \NST_{v_n}(G)} \omega(T)
\sum_{Q\in \D[V(G)\setminus V(T)]}\omega(Q)
\nonumber \\
&=&\sum_{T\in \T(G)}\omega(T)
+\sum_{T\in \NST_{v_n}(G)} \omega(T)
F(G- V(T),\omega).
\end{eqnarray}

Thus Theorem~\ref{identity} is proved.
\proofend

Observe that Theorem~\ref{main}
follows directly from Theorem~\ref{identity}
by taking $u=v_n$ and $y_j=1$ for all $j\in [m]$.

\section{Application}
\label{sec4}

\def \C {{\mathscr C}}

In the last section, we give some examples of applying Theorem~\ref{main}
to determine spanning numbers of graphs.

Let $G$ be a connected multi-graph with $u\in V(G)$.
For $1\le i\le |V(G)|-2$,
let $\C_i(G, u)$ (or simply $\C_i(u)$)
be the set of connected induced subgraphs $G[S]$,
where $u\in S\subset V(G)$, such that $|S|=i$
and $G-S$ has no isolated vertices.
Clearly, $|\C_1(u)|\le 1$ and $|\C_2(u)|\le |N_G(u)|$,
where $N_G(u)$ is the set of neighbors of $u$ in $G$.

\begin{figure}[htp]
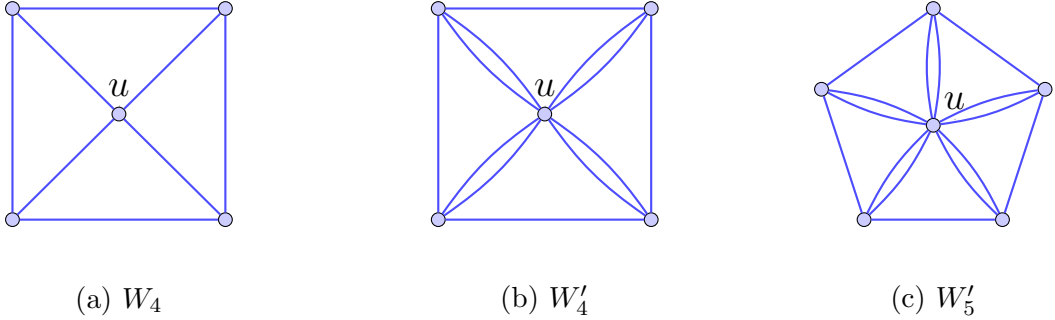

\tikzp{0.7}
{
\foreach \place/\x in {{(-7,0)/1}, {(-7,4)/2},{(-3,4)/3},
{(-3,0)/4}, {(-5,2)/5}}
\node[cblue] (v\x) at \place {};

\filldraw[black] (-5,2.15) circle (0pt) node[anchor=south] {\Large $u$};
\filldraw[black] (-5,-2) circle (0pt) node[anchor=south] {(a) $W_4$};

\path[blue!70,thick] (v5) edge (v1) edge (v2) edge (v3) edge (v4);

\draw[blue!70, thick] (v1) -- (v2) -- (v3)--(v4)--(v1);

\foreach \place/\x in {
{(1,0)/4},{(5,0)/1},
{(3,2)/5},
{(1,4)/3}, {(5,4)/2}}
\node[cblue] (u\x) at \place {};

\filldraw[black] (3,2.15) circle (0pt) node[anchor=south] {\Large $u$};
\filldraw[black] (3,-2) circle (0pt) node[anchor=south] {(b) $W'_4$};


\draw[blue!70, thick] (u1) -- (u2) -- (u3)--(u4)--(u1);

\draw [blue!70, thick, bend left=10, looseness=1] (u5) to (u1);
\draw [blue!70, thick, bend right=10, looseness=1] (u5) to (u1);

\draw [blue!70, thick, bend left=10, looseness=1] (u5) to (u2);
\draw [blue!70, thick, bend right=10, looseness=1] (u5) to (u2);

\draw [blue!70, thick, bend left=10, looseness=1] (u5) to (u3);
\draw [blue!70, thick, bend right=10, looseness=1] (u5) to (u3);

\draw [blue!70, thick, bend left=10, looseness=1] (u5) to (u4);
\draw [blue!70, thick, bend right=10, looseness=1] (u5) to (u4);

\foreach \place/\x in {
{(9,0)/1}, {(11.6,0)/2},
{(10.3,1.788)/3},
{(8.2,2.47)/4},{(12.4,2.47)/5},
{(10.3,4)/6}}
\node[cblue] (b\x) at \place {};

\filldraw[black] (10.7,1.9) circle (0pt) node[anchor=south] {\Large $u$};
\filldraw[black] (10.3,-2) circle (0pt) node[anchor=south] {(c) $W'_5$};


\draw[blue!70, thick] (b1) -- (b2) -- (b5) -- (b6) -- (b4) --(b1);

\draw [blue!70, thick, bend left=10, looseness=1] (b3) to (b1);
\draw [blue!70, thick, bend right=10, looseness=1] (b3) to (b1);
\draw [blue!70, thick, bend left=10, looseness=1] (b3) to (b2);
\draw [blue!70, thick, bend right=10, looseness=1] (b3) to (b2);

\draw [blue!70, thick, bend left=10, looseness=1] (b3) to (b4);
\draw [blue!70, thick, bend right=10, looseness=1] (b3) to (b4);

\draw [blue!70, thick, bend left=10, looseness=1] (b3) to (b5);
\draw [blue!70, thick, bend right=10, looseness=1] (b3) to (b5);

\draw [blue!70, thick, bend left=10, looseness=1] (b3) to (b6);
\draw [blue!70, thick, bend right=10, looseness=1] (b3) to (b6);

}
\caption{Graphs $W_4$, $W'_4$ and $W'_5$}\label{f3-1}
\end{figure}

Observe that expression (\ref{main-eq1}) in
Theorem~\ref{main} is
equivalent to the following one:
\equ{cor}
{
\tau(G)=\prod_{v\in V(G)-\{u\}}d_G(v)
-\sum_{i=1}^{|V(G)|-2}
\sum_{H\in \C_i(u)}
\left (\tau(H)
\prod_{v\in V(G)-V(H)}d_{G-V(H)}(v)\right ).
}
Now we apply (\ref{cor}) to determine
$\tau(W_{4})$, $\tau(W'_{4})$ and $\tau(W'_{5})$,
where $W_4$ is the wheel of order $5$
and $W'_{4}$ and $W'_{5}$ are multi-graphs which can be
obtained  from $W_4$ and $W_5$ respectively by adding new edges
parallel to edges
incident with the central vertex,
as shown in Figure~\ref{f3-1} (b) and (c).

Let $u$ be the central vertex in $W_4$ as shown in
Figure~\ref{f3-1} (a).
By (\ref{cor}), 
we have
\equ{EXM1}
{
\tau(W_4)=3^4-2^4-4\times 1\times (2\times 1^2)
-4\times 3\times 1^2
=45.
}
The above equality follows from the fact that
$|\C_1(u)|=1$, $|\C_2(u)|=|\C_3(u)|=4$,
$\tau(H)=1$ and $G-V(H)\cong C_4$ for $H\in \C_1$,
$\tau(H)=i^{i-2}$
and $G-V(H)$ is a path of length $5-i$
for each $H\in \C_i(u)$ and $i=2,3$.
Again, taking $u$ to be the central vertex in $W'_4$,
we have
\equ{EXM2}
{
\tau(W'_4)=4^4-2^4-4\times 2\times (2\times 1^2)
-4\times 8\times 1^2=192.
}
The above equality follows from the fact that
$|\C_1(u)|=1$, $|\C_2(u)|=|\C_3(u)|=4$,
$\tau(H)=1$ and $G-V(H)\cong C_4$ for $H\in \C_1$,
$\tau(H)=2$ for each $H\in \C_2(u)$,
$\tau(H)=8$ for each $H\in \C_3(u)$,
and $G-V(H)$ is a path of length $5-i$
for each $H\in \C_i(u)$ and $i=2,3$.

Similarly, taking $u$ to be the central vertex in $W'_5$,
we have
\equ{EXM3} 
{
\tau(W'_5)=4^5-2^5-5\times 2\times (2\times 2\times 1^2)
-5\times 8 \times 2
-5(4\times 3^2 - 2 - 2\times2)=722.
}
The above equality follows from the fact that
$|\C_1(u)|=1$, $|\C_i(u)|=5$ for $i=2,3,4$,
$\tau(H)=1$ and $G-V(H)\cong C_5$ for $H\in \C_1$,
$\tau(H)=2$ for each $H\in \C_2(u)$,
$\tau(H)=8$ for each $H\in \C_3(u)$,
$\tau(H)=4\times 3^2 - 2 - 2=32$ for each $H\in \C_4(u)$,
and $G-V(H)$ is a path of length $6-i$
for each $H\in \C_i(u)$ and $i=2,3,4$.

Our examples above show that
as an alternative method of computing spanning trees in small graphs
by hand,
applying Theorem~\ref{main}
is sometimes not less efficient than other methods.

Another potential usefulness of this formula is, maybe for some graph classes,
we can use Theorem \ref{main} to obtain 
a better upper bound for the
number of spanning trees than Theorem \ref{Thomassen}.
Corollary \ref{cor1} below is an example. 

\begin{corollary}\label{cor1}
Let $G$ be a graph with degree sequence $d_1\leq d_2\leq \cdots\leq d_n$. Then 
$$
\tau(G)\leq \prod_{i=1}^{n-1} d_i-\prod_{i=1}^{n-1} (d_i-1).
$$
\end{corollary}

\section*{Acknowledgements}

The authors would like to thank the referees for their
constructive  comments.

\end{document}